\documentclass[12points]{article}

\usepackage{amsmath}
\usepackage{amssymb}

\begin{document}

\title{Derivatives with respect to the order \\ of the Bessel function of the first kind}

\author{J. Sesma \thanks{Email: javier@unizar.es}  \\
{\em Departamento de F\'{\i}sica Te\'orica, Facultad de Ciencias,} \\{\em 50009, Zaragoza, Spain}}

\maketitle

\begin{abstract}
An explicit expression of the $k$-th derivative of the Bessel function $J_\nu(z)$, with respect to its order $\nu$, is given. Particularizations for the cases of positive or negative integer $\nu$ are considered.
\end{abstract}

\bigskip

\section{Introduction}

Along this paper we use the notation
\begin{equation}
\mathcal{G}^{(k)}(t) \equiv \frac{d^k}{dt^k}\,\frac{1}{\Gamma(t)}\,,  \quad
\mathcal{P}_m^{(k)}(t) \equiv \frac{1}{k!}\,\frac{d^k}{dt^k}\,(t)_m\,,    \quad
\mathcal{Q}_{m}^{(k)}(t) \equiv \frac{1}{k!}\,\frac{d^k}{dt^k}\,\frac{1}{(t)_m}\,,  \label{i1}
\end{equation}
to refer to the  derivatives of the reciprocal gamma function and of the Pochhammer and reciprocal Pochhammer symbols.

Our purpose is to provide with a closed expression for the $k$-th derivative of the Bessel function $J_\nu(z)$ with respect to its order $\nu$, that we assume to be real. From the ascending series definition \cite[Eq. 10.2.2]{nist}
\begin{equation}
J_\nu(z)=(z/2)^\nu\,\sum_{m=0}^\infty \frac{(-\,z^2/4)^m}{m!\,\Gamma(\nu\! +\! 1\! +\! m)},  \label{i3}
\end{equation}
one obtains immediately, with the notation introduced in (\ref{i1}),
\begin{equation}
\frac{\partial}{\partial \nu}\,J_\nu(z) = J_\nu(z)\,\ln(z/2)+(z/2)^\nu\,\sum_{m=0}^\infty \mathcal{G}^{(1)}(\nu\! +\! 1\! +\! m)\,
\frac{(-\,z^2/4)^m}{m!},  \label{i4}
\end{equation}
an expression that can be found in all treatises dealing with Bessel functions. (See, for instance, \cite[Eq. 10.15.1]{nist}.)
Derivation, $k$-1 times, with respect to $\nu$ gives a recurrence relation
\begin{eqnarray}
\frac{\partial^k}{\partial \nu^k}\,J_\nu(z) &=& \ln(z/2)\left(\dfrac{\partial^{k-1}}{\partial \nu^{n-1}}\,J_\nu(z)\right)+(z/2)^\nu\,\sum_{m=0}^\infty \Bigg[ \dfrac{(-z^2/4)^m}{m!}  \nonumber  \\
& &\hspace{40pt}\times\,\sum_{l=1}^k {k\! -\! 1 \choose l\! -\! 1}\mathcal{G}^{(l)}(\nu\! +\! 1\! +\! m)\,(\ln(z/2))^{k-l}\Bigg],  \label{i5}
\end{eqnarray}
which would allow one to compute the successive derivatives to get the $k$-th one. We suggest, however, a procedure to obtain directly such $k$-th derivative, without need of computing the lower ones. As auxiliary for the main result, we consider, respectively in Sections 2, 3 and 4, explicit expressions for the symbols defined in Eqs. (\ref{i1}). Then, in Section 5, the $k$-th derivative of $J_{\nu}(z)$ with respect to $\nu$ is discussed. The possibility of extending the resulting expressions to the case of complex $\nu$ is discussed in Section 6.

\section{Derivatives of the reciprocal gamma function}

We start with the series expansion \cite[Eq. 5.7.1]{nist}
\begin{equation}
\frac{1}{\Gamma(t)}=\sum_{j=1}^\infty c_j\,t^j,   \label{ii2}
\end{equation}
convergent in the whole complex $t$-plane. Term by term derivation gives
\begin{equation}
\mathcal{G}^{(k)}(t)=\sum_{j=k}^\infty c_j\,\frac{j!}{(j-k)!}\,t^{j-k},    \label{ii3}
\end{equation}
an expansion also convergent in the whole plane. Nevertheless, its convergence becomes slower and slower as $k$ or $|t|$ increase. It is not recommended, for numerical computation, unless $|t|<1$. For large values of $|t|$ it is preferable to use the asymptotic expansion obtained, in a former paper \cite[Appendix B]{abad}, by application of the saddle point method \cite[Sec. 3.6.3]{temm} to the Hankel contour representation of the reciprocal Gamma function.

The expressions of the derivatives of the Bessel function, to be given below, contain $\mathcal{G}^{(k)}(1+\varepsilon)$, which can be calculated by using
\begin{equation}
\mathcal{G}^{(k)}(1+\varepsilon)=\sum_{j=0}^\infty c_{j+k+1}\,(j+1)_k\,\varepsilon^{j},  \label{ii4}
\end{equation}
whenever $|\varepsilon|<1$.

\section{Derivatives of the Pochhammer symbol $(t)_m$ with respect to its argument}

As $(t)_m$ is a polynomial of degree $m$ in $t$, derivatives of order greater than $m$ vanish,
\begin{equation}
\mathcal{P}_m^{(k)}(t) = 0 \qquad\mbox{for}\qquad k>m\,.  \label{iii2}
\end{equation}
For the nontrivial case of $k\leq m$, the $\mathcal{P}_m^{(k)}(t)$ can be computed by means of the recurrence relation
\begin{equation}
\mathcal{P}_{m+1}^{(k)}(t)=(t+m)\,\mathcal{P}_m^{(k)}(t)+\mathcal{P}_m^{(k-1)}(t), \qquad k>0\,, \label{iii9}
\end{equation}
with starting values
\begin{equation}
\mathcal{P}_{0}^{(k)}(t)=\delta_{k,0}\,, \qquad \mathcal{P}_m^{(0)}(t)=(t)_m\,. \label{iii10}
\end{equation}

Explicit expressions for $\mathcal{P}_m^{(k)}(t)$ can be found easily. From the generating function
of the Pochhammer symbols \cite[Sec. 6.2.1, Eq. (2)]{luke}
\begin{equation}
\sum_{m=0}^\infty \,(t)_m \,(-z)^m/m! \equiv \ _1\!F_0(t;;-z) = (1+z)^{-t}, \qquad |z|<1, \label{iii3}
\end{equation}
one obtains by derivation, $k$ times, with respect to $t$
\begin{equation}
\sum_{m=0}^\infty k!\, \mathcal{P}_m^{(k)}(t)\, (-z)^m/m! = (-1)^k\,(1+z)^{-t}\, \left[\ln(1+z)\right]^k, \qquad |z|<1. \label{iii4}
\end{equation}
The term $\mathcal{P}_m^{(k)}(t)$ can be isolated in this way
\begin{eqnarray}
\hspace{-1cm}\mathcal{P}_m^{(k)}(t) &=& \frac{(-1)^{k-m}}{k!}\,\left.\frac{\partial^m}{\partial z^m}\left((1+z)^{-t}\, \left[\ln(1+z)\right]^k\right)\right|_{z=0}  \nonumber \\
&=& \frac{(-1)^{k-m}}{k!}\,\sum_{l=0}^m {m \choose l}\left.\left(\frac{\partial^l}{\partial z^l}(1+z)^{-t}\right)\left(\frac{d^{m-l}}{d z^{m-l}}\left[\ln(1+z)\right]^k\right)\right|_{z=0}\hspace{-12pt}. \label{iii5}
\end{eqnarray}
Now we make use of the trivial result
\begin{equation}
\left.\frac{\partial^l}{\partial z^l}(1+z)^{-t}\right|_{z=0} = (-1)^l\,(t)_l  \label{iii6}
\end{equation}
and of the  generating relation of the Stirling numbers of the first kind \cite[Eq. 26.8.8]{nist},
\begin{equation}
[\ln (1+z)]^k = k! \sum_{n=k}^\infty s(n, k)\, z^n/n!\,, \qquad |z|<1,  \label{iii7}
\end{equation}
to obtain the explicit expression
\begin{equation}
\mathcal{P}_m^{(k)}(t) = (-1)^{m-k}\sum_{l=0}^{m-k}(-1)^l{m \choose l}\,s(m\!-\!l, k)\,(t)_l \qquad\mbox{for}\quad m\geq k\,. \label{iii8}
\end{equation}
An alternative expression, in terms of generalized Bernoulli polynomials \cite{bry1,bry2,sriv},
\begin{equation}
\mathcal{P}_m^{(k)}(t) = (-1)^{m-k}\,{m \choose k}\,B_{m-k}^{(m+1)}(1\! -\! t)\,.  \label{A19}
\end{equation}
can be obtained from a recent paper by Coffey \cite[Eq. (2.5)]{coff}.

For the particular case of $t=0$, Eq. (\ref{iii8}) gives
\begin{eqnarray}
\mathcal{P}_0^{(k)}(0)&=& \delta_{k,0}\,,  \qquad \mathcal{P}_m^{(0)}(0)= \delta_{m,0}\,,   \label{A5}  \\
\mathcal{P}_m^{(k)}(0)&=& (-1)^{m-k}\,s(m, k) \qquad \mbox{for}\quad m\geq k>0\,. \label{A6}
\end{eqnarray}
In the case of $t=1$, use of the property \cite[Eq. 26.8.20]{nist}
\begin{equation}
s(n\!+\!1,k\!+\!1)=n!\sum_{j=k}^n \frac{(-1)^{n-j}}{j!}\,s(j, k)   \label{A8}
\end{equation}
allows one to obtain, from Eq. (\ref{iii8}),
\begin{equation}
\mathcal{P}_m^{(k)}(1) = (-1)^{m-k}\,s(m\!+\!1, k\!+\!1)\,.  \label{A9}
\end{equation}

\section{Derivatives of the reciprocal Pochhammer symbol $1/(t)_m$ with respect to its argument}

For numerical implementation of the derivatives of the reciprocal Pochhamer symbol with respect to its variable, one may use the recurrence relations
\begin{equation}
\mathcal{Q}_{m+1}^{(k)}(t)=\left(\mathcal{Q}_{m}^{(k)}(t) - \mathcal{Q}_{m+1}^{(k-1)}(t)\right)/(t+m), \label{iv4}
\end{equation}
with initial values
\begin{equation}
\mathcal{Q}_0^{(k)}(t)=\delta_{k,0}\,, \qquad   \mathcal{Q}_{m}^{(0)}(t)=1/(t)_m\,, \label{iv5}
\end{equation}

Very simple explicit expressions of the $\mathcal{Q}_{m}^{(k)}(t)$ can be easily obtained from the relation \cite[Eq. 4.2.2.45]{prud}
\begin{equation}
\frac{1}{(t)_m} = \sum_{l=0}^{m-1}\frac{(-1)^l}{l!\,(m-1-l)!}\,\frac{1}{t+l}\,, \qquad m>0\,,   \label{iv2}
\end{equation}
provided $t$ is different from a nonpositive integer, $-n$, such that $0\leq n<m$.
Direct derivation with respect to $t$ in this equation gives
\begin{equation}
\mathcal{Q}_0^{(k)}(t)=\delta_{k,0}, \qquad \mathcal{Q}_{m}^{(k)}(t)= (-1)^k \sum_{l=0}^{m-1}\frac{(-1)^l}{l!\,(m-1-l)!}\,\frac{1}{(t+l)^{k+1}}\,.  \label{iv3}
\end{equation}
For the particular case of $t=1$, this expression admits a more concise form in terms of {\em modified} generalized harmonic numbers, $\hat{H}_m^{(k)}$, defined by
\begin{equation}
\hat{H}_0^{(k)}\equiv\delta_{k,0}\,,\qquad \hat{H}_m^{(k)} \equiv \sum_{j=1}^m \,(-1)^{j-1}\,{m \choose j}\,\frac{1}{j^k}\,, \qquad m\geq 1\,,   \label{A11}
\end{equation}
not to be confused with the generalized harmonic numbers,
\begin{equation}
H_m^{(k)} \equiv \sum_{j=1}^m \frac{1}{j^k}\,, \qquad m\geq 1\,,  \label{A10}
\end{equation}
although
\begin{equation}
\hat{H}_m^{(1)} = {H}_m^{(1)} \equiv H_m  \qquad {\rm for}\quad m\geq 1\,.  \label{H1}
\end{equation}
Besides the explicit expression (\ref{A11}), the recurrence relation
\begin{equation}
\hat{H}_{m+1}^{(k)}=\hat{H}_{m}^{(k)}+\frac{1}{m\!+\!1}\, \hat{H}_{m}^{(k-1)}\,, \qquad m\geq 0\,, \quad k\geq 1\,,  \label{H2}
\end{equation}
with the starting values
\begin{equation}
\hat{H}_{0}^{(k)}=\delta_{k,0}\,, \qquad   \hat{H}_{m}^{(0)}=1\,,  \label{H3}
\end{equation}
may be used to calculate the $\hat{H}_m^{(k)}$. With that notation, Eq. (\ref{iv3}) gives
\begin{equation}
\mathcal{Q}_m^{(k)}(1)=\frac{(-1)^k}{m!}\,\hat{H}_m^{(k)}\,.  \label{A12}
\end{equation}

\section{Derivatives of $J_\nu(z)$ with respect to $\nu$}

We proceed to obtain our expression for the $k$-th derivative of $J_\nu(z)$ with respect to $\nu$. To avoid unnecessary complications in the resulting formulas, we assume $k\neq 0$, i. e., $k= 1, 2, \ldots$.

Let us denote by $N$ the nearest integer to $\nu$, and define $\varepsilon$ by
\begin{equation}
\nu=N+\varepsilon, \qquad  |\varepsilon|\leq 1/2.  \label{v1}
\end{equation}
We distinguish two possible ranges of values of $N$.

\subsection{$N\geq 0$}

The ascending series in Eq. (\ref{i3}) can be written in the form
\begin{equation}
J_\nu(z)=(z/2)^\nu\,\frac{1}{\Gamma(1+\varepsilon)}\sum_{m=0}^\infty \frac{(-\,z^2/4)^m}{m!\,(1+\varepsilon)_{m+N}}\,.  \label{v2}
\end{equation}
Derivation, $k$ times, with respect to $\nu$ gives, with the notation introduced in (\ref{i1}),
\begin{eqnarray}
\frac{\partial^k}{\partial \nu^k}\,J_\nu(z) &=& k!\,(z/2)^\nu\,\sum_{m=0}^\infty  \frac{(-z^2/4)^m}{m!}  \nonumber  \\
& & \times\,\sum_{k_1=0}^k \frac{\left[\ln(z/2)\right]^{k_1}}{k_1!}\sum_{k_2=0}^{k-k_1} \frac{\mathcal{G}^{(k_2)}(1\!+\!\varepsilon)}{k_2!}\,\mathcal{Q}_{m+N}^{(k-k_1-k_2)}(1\!+\!\varepsilon)\,,  \label{v3}
\end{eqnarray}
where $\mathcal{G}^{(k_2)}(1\!+\!\varepsilon)$ is given in Eq. (\ref{ii4}) and, according to Eq (\ref{iv3}),
\begin{equation}
\mathcal{Q}_0^{(k)}(1\!+\!\varepsilon)= \delta_{k,0}\,,\quad \mathcal{Q}_{m+N}^{(k)}(1\!+\!\varepsilon)=\sum_{j=1}^{m+N}\frac{(-1)^{k+j-1}}{(j\!-\!1)!\,(m\!+\!N\!-\!j)!}\,\frac{1}{(\varepsilon\!+\!j)^{k+1}}\,. \label{v4}
\end{equation}
In the particular case of $\nu$ being a nonnegative integer, $\nu=n\geq 0$, Eq. (\ref{v3}) becomes, in terms of the modified generalized harmonic numbers defined in (\ref{A11}),
\begin{eqnarray}
\left.\frac{\partial^k}{\partial \nu^k}\,J_\nu(z)\right|_{\nu=n} &=& k!\,(z/2)^n\,\sum_{m=0}^\infty  \frac{(-z^2/4)^m}{m!\,(m\!+\!n)!} \nonumber  \\
&&  \hspace{-20pt}\times\sum_{k_1=0}^k \frac{\left[\ln(z/2)\right]^{k_1}}{k_1!}\sum_{k_2=0}^{k-k_1} (-1)^{k-k_1-k_2}\,c_{k_2+1}\,\hat{H}_{m+n}^{(k-k_1-k_2)}.  \label{v5}
\end{eqnarray}

Expressions for the first derivative can be found in the bibliography. Besides the familiar expressions given in, for instance, Sect. 10.15 of Ref. \cite{nist}, alternative closed forms can be found in a paper by Brychkov and Geddes  \cite{bry3}. Our Eqs. (\ref{v3}) and (\ref{v5}) become, for $k=1$,
\begin{eqnarray}
\frac{\partial}{\partial \nu}\,J_\nu(z) &=& \left(\ln (z/2)-\psi(1\!+\!\varepsilon)\right)\,J_\nu(z)   \nonumber  \\
&& \hspace{-20pt}+\,\frac{(z/2)^\nu}{\Gamma(1\!+\!\varepsilon)}\,\sum_{m=0}^\infty  \frac{(-z^2/4)^m}{m!}
 \sum_{j=1}^{m+N}\frac{(-1)^{j}}{(j\!-\!1)!\,(m\!+\!N\!-\!j)!}\,\frac{1}{(\varepsilon\!+\!j)^2},  \label{v7}
\end{eqnarray}
where $\psi$ represents the digamma function and the last sum is understood to be zero if $m+N=0$. In the case of integer $\nu=n\geq 0$ we have
\begin{equation}
\left.\frac{\partial}{\partial \nu}\,J_\nu(z)\right|_{\nu=n} = \left(\ln (z/2)+\gamma\right)\,J_n(z) - (z/2)^n
\sum_{m=0}^\infty  \frac{(-z^2/4)^m}{m!\,(m\!+\!n)!}\,\hat{H}_{m+n}^{(1)}\,,  \label{v8}
\end{equation}
where $\gamma$ represents the well known Euler-Mascheroni constant.

\subsection{$N<0$}

Instead of Eq. (\ref{v2}) we have now
\begin{eqnarray}
J_\nu(z)&=&(z/2)^\nu\,\frac{1}{\Gamma(1+\varepsilon)}\Bigg[\sum_{m=0}^{-N-1} \frac{(-\,z^2/4)^m}{m!}\,(-1)^{-N-m}\,(-\varepsilon)_{-N-m}  \nonumber \\
& & \hspace{100pt}+\, \sum_{m=-N}^\infty \frac{(-\,z^2/4)^m}{m!}\,\frac{1}{(1+\varepsilon)_{m+N}}\Bigg].   \label{v9}
\end{eqnarray}
Derivation with respect to $\nu$ gives
\begin{eqnarray}
\frac{\partial^k}{\partial\nu^k}J_\nu(z)&=&k!\,(z/2)^\nu\sum_{k_1=0}^k \frac{\left[\ln (z/2)\right]^{k_1}}{k_1!}
\sum_{k_2=0}^{k-k_1}\frac{\mathcal{G}^{(k_2)}(1+\varepsilon)}{k_2!}    \nonumber \\
& & \hspace{10pt} \times\,\Bigg[\sum_{m=0}^{-N-1} \frac{(-z^2/4)^m}{m!}\,(-1)^{-N-m+k-k_1-k_2}\,\mathcal{P}_{-N-m}^{(k-k_1-k_2)}(-\varepsilon) \nonumber   \\
& & \hspace{60pt}+\, \sum_{m=-N}^\infty \frac{(-\,z^2/4)^m}{m!}\,\mathcal{Q}_{m+N}^{(k-k_1-k_2)}(1+\varepsilon)\Bigg],  \label{v10}
\end{eqnarray}
with $\mathcal{P}_{-N-m}^{(k)}(-\varepsilon)$ given by Eqs. (\ref{iii8}) or (\ref{A19}) and $\mathcal{Q}_{m+N}^{(k)}(1+\varepsilon)$ by Eq. (\ref{v4}). In the particular case of $\nu$ being a negative integer, $\nu=-n$, $n>0$, this equation turns into
\begin{eqnarray}
\left.\frac{\partial^k}{\partial\nu^k}J_\nu(z)\right|_{\nu=-n}&=&k!\,(z/2)^{-n}\sum_{k_1=0}^k \frac{\left[\ln (z/2)\right]^{k_1}}{k_1!}
\sum_{k_2=0}^{k-k_1}\,c_{k_2+1}    \nonumber \\
& & \times\,\Bigg[\sum_{m=0}^{n-1} \frac{(-z^2/4)^m}{m!}\,s(n\!-\!m,\,k\!-\!k_1\!-\!k_2) \nonumber   \\
& & \hspace{10pt}+\, (-1)^{k-k_1-k_2}\,\sum_{m=n}^\infty \frac{(-\,z^2/4)^m}{m!\,(m-n)!}\,\hat{H}_{m-n}^{(k-k_1-k_2)}\Bigg].  \label{v11}
\end{eqnarray}
For $k=1$, Eq. (\ref{v10}) becomes
\begin{eqnarray}
\frac{\partial}{\partial \nu}\,J_\nu(z) &=& \left(\ln (z/2)-\psi(1+\varepsilon)\right)\,J_\nu(z)   \nonumber  \\
& & \hspace{-10pt}+\,\frac{(z/2)^\nu}{\Gamma(1\!+\!\varepsilon)}
\Bigg[\sum_{m=0}^{-N-1} \frac{(-z^2/4)^m}{m!}\,(-N\!-\!m)\,B_{-N-m-1}^{(-N-m+1)}(1\!+\!\varepsilon)   \nonumber  \\
& & +\sum_{m=-N}^\infty  \frac{(-z^2/4)^m}{m!}\,
 \sum_{j=1}^{m+N}\frac{(-1)^{j}}{(j\!-\!1)!\,(m\!+\!N\!-\!j)!}\,\frac{1}{(\varepsilon\!+\!j)^2}\Bigg],  \label{v12}
\end{eqnarray}
where $B_n^{(\alpha)}(x)$ represents the generalized Bernoulli polynomial \cite{bry1,bry2,sriv}. It may be written in terms of Stirling numbers of the first kind by using the relation
\begin{equation}
(-N\!-\!m)\,B_{-N-m-1}^{(-N-m+1)}(1\!+\!\varepsilon) = \sum_{j=0}^{-N-m-1}\,(j\!+\!1)\,s(-N\!-\!m, j\!+\!1)\,\varepsilon^j\,.  \label{v13}
\end{equation}
In the case of $\nu$ being a negative integer, Eq. (\ref{v12}) gives
\begin{eqnarray}
\left.\frac{\partial}{\partial\nu}J_{\nu(z)}\right|_{\nu=-n} &=&\left(\ln (z/2)+\gamma\right)\,J_{-n}(z)  \nonumber  \\
&&\hspace{-90pt}-\,(z/2)^{-n}\Bigg[(-1)^{n}\sum_{m=0}^{n-1} \frac{(z^2/4)^m}{m!}\,(n\!-\!m\!-\!1)!  +  \sum_{m=n}^\infty \frac{(-\,z^2/4)^m}{m!\,(m-n)!}\,\hat{H}_{m-n}^{(1)}\Bigg].   \label{v14}
\end{eqnarray}

\section{Extension to complex values of $\nu$}

The expressions of the derivatives of the reciprocal Gamma function and of the Pochhammer and reciprocal Pochhammer symbols given in sections 2 to 4
stay for complex values of their argument $t$. Therefore, our Eqs. (\ref{v3}), (\ref{v7}), (\ref{v10}) and (\ref{v12}) may be used safely for complex $\varepsilon$, i. e. complex $\nu$, whenever $|\Im \nu|\lesssim 1/2$. As auxiliary integer $N$ one should consider again the nearest to $\nu$ one, in such a way that, instead of (\ref{v1}), one would have
\begin{equation}
\nu=N+\varepsilon, \qquad  |\Re \varepsilon|\leq 1/2.  \label{vi1}
\end{equation}
For large values of $\Im \nu$, the given expressions are correct, but they are not useful from a computational point of view. The reason, as pointed out is Sect. 2, is the slow convergence of the series in the right hand side of (\ref{ii3}) for large values ot $t$.

\section*{Acknowledgements}

This work has been supported by Departamento de Ciencia, Tecnolog\'{\i}a y Universidad del Gobierno de Arag\'on (Project E24/1) and Ministerio de Ciencia e Innovaci\'on (Project MTM2009-11154)

\end{document}